\documentstyle[11pt]{article}

\newtheorem{eg}{{\bf Example}}
\newtheorem{lemma}{{\bf Lemma}}[section]
\newtheorem{theo}{{\bf Theorem}}
\newtheorem{prop}{{\bf Proposition}}

\newtheorem{defn}{{\bf Definition}}
\newtheorem{remark}{{\bf Remark}}
\newtheorem{conj}{{\bf Conjecture}}

\font\bbb=msbm10 scaled\magstep1

\newcommand{\CC}{\mbox{\bbb C}}

\newcommand{\RR}{\mbox{\bbb R}}
\newcommand{\ZZ}{\mbox{\bbb Z}}

\def\scoll{\mbox{{\scriptsize $\searrow$}$\!\!\stackrel{{}^{\rm s}}{}\,\,$}}

\def\coll{\mbox{\scriptsize $\,\searrow\,$}}

\newcommand{\La}{\Lambda}

\newcommand{\si}{\sigma}

\begin{document}

\title{\bf Non-existence of 6-dimensional pseudomanifolds with complementarity}
\author{{\bf Bhaskar Bagchi} and {\bf Basudeb Datta}
} \vspace{-5mm}
\date{To appear in Advances in Geometry}

\maketitle

{\small

\noindent {\bf Abstract.} In a previous paper (\cite{d3}) the
second author showed that if $M$ is a  pseudomanifold with
complementarity other than the 6-vertex real projective plane and
the 9-vertex complex projective plane, then $M$ must have
dimension $\geq 6$, and - in case of equality - $M$ must have
exactly 12 vertices. In this paper we prove that such a
6-dimensional pseudomanifold does not exist. On the way to proving
our main result we also prove that all combinatorial
triangulations of the 4-sphere with at most 10 vertices are
combinatorial $4$-spheres. }

\bigskip

{\small


\smallskip

\mbox{}
\hspace{16mm} complementarity, piecewise-linear manifolds.

}

\section{Introduction and results}

Recall that a {\em simplicial complex} is a collection of
non-empty finite sets such that every non-empty subset of an
element is also an element.  For $i \geq 0$, the elements of size
$i+1$ are called the {\em $i$-simplices} or {\em $i$-faces} of the
complex. The vertex-set $V(K)$ of $K$ is by definition the union
of all the faces in $K$. We identify the 0-faces with the
vertices. All the simplicial complexes considered in this paper
are with finite vertex-set. For a simplicial complex $K$, the
maximum of $k$ such that $K$ has a $k$-face is called the {\em
dimension} of $K$. An 1-dimensional simplicial complex is also
called a {\em graph}.

If $K$, $L$ are two simplicial complexes, then a {\em simplicial
isomorphism} from $K$ to $L$ is a bijection $\pi : V(K) \to V(L)$
such that for $\si\subseteq V(K)$, $\si$ is a face of $K$ if and
only if $\pi (\si)$ is a face of $L$. The complexes $K$, $L$ are
called (simplicially) {\em isomorphic} when such an isomorphism
exists. We identify two simplicial complexes if they are
isomorphic.

A simplicial complex $K$ is usually thought of as a prescription
for constructing a topological space (polyhedron), called the {\em
geometric carrier} of $K$ and denoted by $|K|$, by pasting
together geometric simplices. We say that a simplicial complex $K$
{\em triangulates} a topological space $X$ (or $K$ is a {\em
triangulation} of $X$) if $X$ is homeomorphic to $|K|$. A
simplicial complex $K$ is said to be {\em connected}
(respectively, {\em simply connected}) if $|K|$ is connected
(respectively, simply connected). Similarly, $K$ is said to be
{\em contractible} if $|K|$ is contractible.

If $\si$ is a face of a simplicial complex $K$ then the {\em link}
of $\si$ in $K$, denoted by ${\rm Lk}_K(\sigma)$ (or simply by
${\rm Lk}(\sigma)$), is by definition the complex whose faces are
the faces $\tau$ of $K$ such that $\tau$ is disjoint from $\si$
and  $\si\cup\tau$ is a face of $K$.

If the number of $i$-simplices of a $d$-dimensional simplicial
complex $K$ is $f_i(K)$ ($0\leq i\leq d$), then the number
$\chi(K) := \sum_{i=0}^{d}(-1)^if_i(K)$ is called the {\em Euler
characteristic}  of $K$. A simplicial complex $K$ is called {\em
$k$-neighbourly} if $f_{k-1}(K)= {f_0(K) \choose k}$.

A simplicial complex $K$ is called {\em pure} if all the maximal
faces of $K$ have the same dimension. A maximal face in a pure
simplicial complex is also called a {\em facet}. With each
$d$-dimensional pure simplicial complex $K$ we associate a graph
$\La (K)$ whose vertices are the facets of $K$ and two such
vertices are adjacent in $\La (K)$ if and only if the
corresponding facets intersect in a $(d-1)$-face. If $\La (K)$ is
connected then $K$ is called {\em strongly connected}.

A $d$-dimensional pure simplicial complex $K$ is called a {\em
weak pseudomanifold} if each $(d-1)$-face is contained in exactly
two $d$-faces of $K$. A strongly connected weak pseudomanifold is
called a {\em pseudomanifold}. A $d$-dimensional pure simplicial
complex $K$ is called a {\em weak pseudomanifold with boundary} if
each $(d-1)$-face is contained in one or two $d$-faces of $K$ and
there exists a $(d-1)$-face of $K$ which is in only one $d$-face
of $K$.

By a {\em subdivision} of a simplicial complex $K$ we mean a
simplicial complex $K^{\prime}$ together with a homeomorphism from
$|K^{\prime}|$ onto $|K|$ which is facewise linear. Two simplicial
complexes $K$ and $L$ are called {\em combinatorially equivalent}
(denoted by {\em $K \approx L$}) if they have isomorphic
subdivisions. So, $K\approx L$ if and only if $|K|$ and $|L|$ are
pl homeomorphic (\cite{rs}).

For a set $V$ with $d + 2$ elements, let $S$ be the simplicial
complex whose faces are all the non-empty proper subsets of $V$.
Then $S$ triangulates the $d$-sphere. This complex is called the
{\em standard $d$-sphere} and is denoted by $S^{\,d}_{d + 2}(V)$
(or simply by $S^{\,d}_{d + 2}$). A polyhedron is called a {\em pl
$d$-sphere} if it is pl homeomorphic to $|S^{\,d}_{d + 2}|$. A
simplicial complex  $X$ is called a {\em combinatorial $d$-sphere}
if it is combinatorially equivalent to $S^{\,d}_{d + 2}$.  So, $X$
is a combinatorial $d$-sphere if and only if $|X|$ is a pl
$d$-sphere.

A simplicial complex $K$ is called a {\em combinatorial
$d$-manifold} if the link of each vertex is a combinatorial
$(d-1)$-sphere. So, $K$ is a combinatorial $d$-manifold if and
only if $|K|$ is a closed pl $d$-manifold (\cite{rs}). If $K$
triangulates $X$ and $K$ is a combinatorial manifold then $K$ is
called a {\em combinatorial triangulation} of $X$. A combinatorial
manifold is automatically a weak pseudomanifold. Further, a
connected combinatorial manifold is a pseudomanifold.

If $M$ is a triangulation of a 2-manifold then the link of a
vertex is a circle and hence $M$ is combinatorial manifold. Again,
the link of a vertex in a triangulation of a 3-manifold is a
triangulation of the 2-sphere and all triangulations of the
2-sphere are combinatorial 2-spheres. So, any triangulation of a
3-manifold is a combinatorial triangulation.

Eells and Kuiper defined a {\em manifold like a projective plane}
to be a cohomology projective plane over reals, complex numbers,
quaternions or Cayley numbers. It is well known that the
projective planes over reals and complex numbers are the only
manifolds like projective planes of dimensions 2 and 4,
respectively. In \cite{bk1}, Brehm and K\"{u}hnel proved\,:

\begin{prop}$\!\!${\bf .}  \label{brku1}
Let $M$ be an $n$-vertex combinatorial $d$-manifold $(d>0)$.
\begin{enumerate}
\vspace*{-1mm} \item[{\rm (a)}] If $n<3\lceil d/2 \rceil+3$ then
$M$ is a combinatorial $d$-sphere. \vspace*{-1.25mm} \item[{\rm
(b)}] If $n=3d/2 +3$ and $M$ is not a combinatorial $d$-sphere
then $d =  2, 4, 8$ or $16$ and $|M|$ must be a `manifold like a
projective plane'.
\end{enumerate}
\end{prop}

Existence or otherwise of more than one smooth structures on
$S^{\,4}$ is a long standing open problem. In view of the
isomorphism between the group of smooth $4$-spheres and the group
of pl $4$-spheres (see \cite[page 201]{ks}), this problem is
equivalent to the existence problem of a combinatorial
triangulation of $S^{\,4}$ which is not a combinatorial sphere. By
Proposition \ref{brku1}, such a triangulation of $S^{\,4}$ must
have at least ten vertices. Here we prove\,:

\begin{theo}$\!\!${\bf .} \label{t1}
If $M$ is a $10$-vertex combinatorial triangulation of $S^{\,4}$
then $M$ is a combinatorial $4$-sphere.
\end{theo}

This result is an immediate consequence of Lemma \ref{l3.4} which
we need to prove our main result (Theorem \ref{t2}).

\begin{defn}$\!\!${\bf .}  {\rm
A simplicial complex $K$ is said to satisfy {\em complementarity}
if exactly one of each complementary pair of non-empty subsets of
$V(K)$ is a face of $K$. A simplicial complex satisfying
complementarity is said to be a {\em complementary} simplicial
complex. }
\end{defn}

In \cite{am}, Arnoux and Marin  proved the following\,:

\begin{prop}$\!\!${\bf .}  \label{arma}
If $M$ is a combinatorial manifold as in part {\rm (b)} of
Proposition $\ref{brku1}$ then $M$ satisfies complementarity.
\end{prop}

In \cite{d2}, the second author proved  the following converse\,:

\begin{prop}$\!\!${\bf .}  \label{da1}
If an $n$-vertex combinatorial $d$-manifold $(d>0)$ satisfies
complementarity then $d=2$, $4$, $8$ or $16$ and $n=3d/2+3$.
\end{prop}

It is easy to see that complementary simplicial complexes are
plentiful. Indeed, let ${\cal A}$ be a family of non-empty subsets
of a finite set $X$ such that (a) $A_1$, $A_2\in {\cal
A}\Rightarrow A_1\cup A_2\neq X$, and (b)  ${\cal A}$ is maximal
(w.r.t. set-inclusion) subject to (a). Then ${\cal A}$ is a
complementary simplicial complex. In contrast, complementary
(weak) pseudomanifolds are very hard to come by.

\begin{eg}$\!\!${\bf .} \label{e2}
{\rm Here is a $7$-vertex $3$-dimensional complementary weak
pseudomanifold with boundary. Its vertex set is $\ZZ_7$. The
facets are $\{i, i+3, i+5, i+6\}$, $i\in\ZZ_7$. Note that in this
example each edge is in exactly two facets and any two facets have
exactly one edge in common. }
\end{eg}

It is well known (see \cite{bd2, k2}) that there is a unique
6-vertex combinatorial triangulation of $\RR P^{\,2}$. It is also
known (see \cite{bd1, bd3}) that there is a unique 9-vertex
combinatorial triangulation of $\CC P^{\,2}$. In \cite{bk2}, Brehm
and K\"{u}hnel constructed three 15-vertex combinatorial
8-manifolds which are not combinatorial spheres. By Proposition
\ref{arma}, all these combinatorial manifolds satisfy
complementarity. In fact, all the known complementary weak
pseudomanifolds are combinatorial manifolds. In \cite{d3}, the
second  author proved the following\,:

\begin{prop}$\!\!${\bf .} \label{da2}
Let $M$ be an $n$-vertex $d$-dimensional $(d>0)$ complementary
pseudomanifold. If $n\leq d+6$ or $d\leq 6$ then $d= 2, 4$ or $6$
and $n=3{d}/{2} +3$. Moreover, if $d\leq 5$ then $M$ is a
combinatorial manifold.
\end{prop}

The proofs in \cite{d3} show that Proposition \ref{da2} remains
true even if we assume $M$ is a complementary weak pseudomanifold.
In view of this result, we make\,:

\begin{conj}$\!\!${\bf .}
If an $n$-vertex $d$-dimensional complementary weak pseudomanifold
exists then $d$ is even and $n=3d/2+3$.
\end{conj}

\begin{conj}$\!\!${\bf .}
If $M$ is a $d$-dimensional complementary weak pseudomanifold
then its Euler characteristic $\chi(M)$ is given by\,:
$$
\chi(M) = \left\{
     \begin{array}{ll}
    1 & \mbox{if $d\equiv 2$ $($mod $4)$} \\
   3 & \mbox{if $d\equiv 0$ $($mod $4)$}.
  \end{array}
\right.
$$
\end{conj}

\begin{remark}$\!\!${\bf .} \label{re1}
{\rm If $M$ is any $n$-vertex complementary simplicial complex
then the total number $2^{\,n - 1} - 1$ of its faces is odd. Hence
$\chi(M)$ is odd. If, further, $n$ is even then it has $2^{\,n -
2}$ even dimensional faces and $2^{\,n - 2} - 1$ odd dimensional
faces, and hence $\chi(M)=1$ in this case. }
\end{remark}

\begin{remark}$\!\!${\bf .} \label{re2}
{\rm The referee suggested the following short proof of
Proposition \ref{da1}. Let $M$ be an $n$-vertex complementary
combinatorial $d$-manifold. By  Remark \ref{re1}, $|M|$ is not a
sphere. Since $M$ has no $(d+1)$-simplex, complementarity implies
that $M$ is $(n-d-2)$-neighbourly and hence $|M|$ is
$(n-d-4)$-connected. This implies that $n-d-4\leq (d-1)/2$.
Proposition 3 now follows from Proposition 1.}
\end{remark}

In view of Proposition \ref{da2}, $d=6$, $n=12$ are the smallest
parameters for which the existence of an $n$-vertex
$d$-dimensional complementary pseudomanifold was an open problem.
In this article we prove\,:

\begin{theo}$\!\!${\bf .} \label{t2}
There does not exist a $12$-vertex $6$-dimensional weak
pseudomanifold with complementarity.
\end{theo}


\section{Preliminaries and Definitions.}

Let $K$, $L$ be two simplicial complexes with disjoint vertex
sets. Then their {\em join} $K \ast L$ is the simplicial complex
whose faces are those of $K$ and of $L$, and the  unions of faces
of $K$ with faces of $L$. Clearly, if $K$ and $L$ are (weak)
pseudomanifolds, then so is $K\ast L$.

For $n\geq 3$, the combinatorial 1-sphere (circle) with $n$
vertices is the unique 1-dimensional $n$-vertex pseudomanifold and
is denoted by $S^{\,1}_n$.

If $\si$ is an $i$-face in a $d$-dimensional weak pseudomanifold
$K$ then ${\rm Lk}_K(\sigma)$ is a $(d-i-1)$-dimensional weak
pseudomanifold. The number of vertices in ${\rm Lk}(\sigma)$ is
called the {\em degree} of $\sigma$ and is denoted by
$\deg_K(\sigma)$ (or simply by $\deg(\sigma)$).

We need the next two results later.

\begin{prop}  \label{brku2} \mbox{\rm (Brehm and K\"{u}hnel
\cite{bk1})}{\bf .} Let $M$ be an $n$-vertex combinatorial
$d$-manifold $(d>0)$. If $n\leq 2d+3-i$ for some $i < d/2$ then
$|M|$ is $i$-connected in the sense of homotopy.
\end{prop}

It follows from the next example that this result is best possible
for $i=0$.

\begin{eg}$\!\!${\bf .} \label{e1}
{\rm For $d\geq 2$, let $K^{\,d}_{2d+3}$ be the $d$-dimensional
simplicial complex whose vertex set is the vertex set of the
circle $S^{\,1}_{2d+3}$ and the facets are the sets of $d+1$
vertices obtained by deleting an interior vertex from the
$(d+2)$-paths in the circle. The simplicial complex
$K^{\,d}_{2d+3}$ is a combinatorial $d$-manifold (see \cite{k1,
k2}) and $|K^{\,d}_{2d+3}|$ is not simply connected. The space
$|K^{\,3}_9|$ is known as the $3$-dimensional Klein bottle. }
\end{eg}

\begin{prop} \label{alst} \mbox{\rm (Altshuler and Steinberg
\cite{as})}{\bf .} If $M$ is a $9$-vertex combinatorial
$3$-manifold then either $M$ is a combinatorial $3$-sphere or $M$
is isomorphic to $K^{\,3}_9$.
\end{prop}

A subcomplex $L$ of a simplicial complex $K$ is called an {\em
induced} (or {\em full}\,) subcomplex of $K$ if $\sigma\in K$ and
the vertices of $\sigma$ are in $L$ imply $\sigma \in L$.

Let $L\subseteq K$ be simplicial complexes. The {\em simplicial
neighbourhood} of $L$ in $K$ is the subcomplex $N(L, K)$ of $K$
whose maximal simplices are those maximal simplices of $K$ which
intersect $V(L)$. Clearly, $N(L, K)$ is the smallest subcomplex of
$K$ such that its geometric carrier is a topological neighbourhood
of $|L|$ in $|K|$. The induced subcomplex $C(L, K)$ on the
vertex-set $V(K)\setminus V(L)$ is called the {\em simplicial
complement} of $L$ in $K$.

Suppose $P^{\,\prime}\subseteq P$ are polyhedra and $P =
P^{\,\prime} \cup B$, where $B$ is a   pl $(k+1)$-ball. If
$P^{\,\prime} \cap B$ is a   pl $k$-ball then we say that there is
an {\em elementary collapse} of $P$ on $P^{\,\prime}$. We say that
$P$ collapses on $Q$ and write $P\coll Q$ if there exists a
sequence $P = P_0$, $P_1, \dots, P_n = Q$ of polyhedra such that
there is an elementary collapse of $P_{i-1}$ on $P_{i}$ for $1\leq
i\leq n$ (\cite{rs}).

A {\em regular neighbourhood} of a polyhedron $Q$ in a pl
$d$-manifold $M$ is a $d$-dimensional submanifold $N$ with
boundary such that $N \coll Q$ and $N$ is a neighbourhood of $Q$
in $M$.

Let $\tau\subset\sigma$ be two faces of a simplicial complex $K$.
We say that $\tau$ is a {\em free face} of $\sigma$ if $\sigma$ is
the only face of $K$ which properly contains $\tau$. (It follows
that $\dim(\sigma)-\dim(\tau)=1$ and $\sigma$ is a maximal simplex
in $K$.) If $\tau$ is a free face of $\sigma$ then $K^{\,\prime}
:= K \setminus \{\tau, \sigma\}$ is a simplicial complex. We say
that there is an {\em elementary collapse} of $K$ on
$K^{\,\prime}$. We say $K$ {\em collapses} on $L$ and write
$K\scoll L$ if there exists a sequence $K=K_0$, $K_1, \dots$,
$K_n=L$ of simplicial complexes such that there is an elementary
collapse of $K_{i-1}$ on $K_{i}$ for $1\leq i\leq n$ (see
\cite{b}). If $L$ consists of a 0-simplex (a point) we say that
$K$ is {\em collapsible} and write $K \scoll 0$. Clearly, if
$K\scoll L$ then $|K|\coll |L|$ as polyhedra and hence $|K|$ and
$|L|$ have the same homotopy type (see \cite{rs}). So, if a
simplicial complex $K$ is collapsible then $|K|$ is contractible.

\section{Ten-vertex four-sphere.}

\begin{lemma}$\!\!${\bf .} \label{l3.1}
Let $S$ be a combinatorial triangulation of a simply connected
$4$-manifold with $\chi(S)=2$ and $\sigma$ be a facet of $S$. Let
$L$ be the induced subcomplex of $S$ on $V(S)\setminus\sigma$.
Then
\begin{enumerate}
     \item[$(a)$] $L$ is contractible.
     \item[$(b$)] If, further, $L$ is collapsible then $S$ is a
     combinatorial sphere.
     \end{enumerate}
\end{lemma}

\noindent {\bf Proof\,.} Being simply connected, $|S|$ is
orientable. Using Poincar\'{e} duality, one sees that any simply
connected 4-manifold of Euler characteristic 2 is a homotopy
sphere. So, $|S|$ is a homotopy 4-sphere.

Let $D = |S| \setminus |\sigma|$. Then, $D$ is homotopic to a
point and hence contractible.

If $x\in D \setminus |L|$ then there exists a unique pair
$(\alpha, \beta)$, where $\alpha\subseteq\sigma$ and $\beta\in L$
such that $x$ is in the interior of $|\alpha\cup\beta|$. So, there
exist a unique pair $(y, z)\in |\alpha| \times |\beta|$ and
$0<s<1$ such that $x = sy + (1-s)z$. Then $H\colon D \times [0,
1]\to D$, given by
$$
H(x, t) = \left\{
     \begin{array}{ll}
    x & \mbox{if $x \in |L|$} \\
   (1 - t)sy + (1 - s + ts)z & \mbox{if $x = sy + (1-s)z\not\in |L|$},
  \end{array}
\right.
$$
defines a homotopy between $D$ and $|L|$. So, $|L|$ is
contractible. This proves $(a)$.

Let $J = S^{\,3}_5(\sigma)$ and $K = S\setminus\{\sigma\}$. Let $X
= |L|$ and $N=|K|$.

Then (i) $N$ is a neighbourhood of $X$ in the closed pl manifold
$|S|$, (ii) $N$ is a compact pl manifold with boundary (see
\cite[Corollary 3.14]{rs}) and (iii) $(K, L, J)$ is a
triangulation of $(N, X, \partial N)$ with $L$ a full subcomplex
of $K$, $K=N(L, K)$ and $J= N(L, K)\cap C(L, K)$. Therefore by the
Simplicial Neighbourhood Theorem (see \cite[page 34]{rs}), $N$ is
a regular neighbourhood of $X$.

Now, if $L$ is collapsible, then $N$ is a regular neighbourhood of
the collapsible polyhedron $X$ and hence (see \cite[page 41]{rs})
$N$ is a pl 4-ball.

Let $\sigma = 12345$. Then $B := |S^{\,4}_6(\{1, \dots, 6\})
\setminus \{\sigma\}|$ is a pl 4-ball. Let $\varphi \colon N\to B$
be a pl-homeomorphism. Let $u$ be a point in the interior of
$|\sigma|$. For $x\in |\sigma|$ and $x\neq u$ there exist a unique
point $y\in |\sigma|\cap N$ and $0\leq t<1$ such that $x = tu + (1
- t)y$. Then $\widetilde{\varphi}$, given by
$$
\widetilde{\varphi}(x) = \left\{
     \begin{array}{ll}
     u & \mbox{if $x = u$} \\
    tu + (1 - t)\varphi(y) & \mbox{if $x \in |\sigma|$, $x\neq u$
    and $x = tu + (1 - t)y$} \\
   \varphi(x) & \mbox{otherwise},
  \end{array}
\right.
$$
defines a pl-homeomorphism between $|S|$ and $|S^{\,4}_6(\{1,
\dots, 6\})|$. This proves $(b)$. \hfill $\Box$

\begin{lemma}$\!\!${\bf .} \label{l3.2}
Let $N$ be a $2$-dimensional simplicial complex on at most five
vertices. Suppose each edge of $N$ is contained in at least two
triangles of $N$. Then $N$ contains a combinatorial $2$-sphere as
a subcomplex.
\end{lemma}

\noindent {\bf Proof\,.} Clearly, $\#(V(N)) \geq 4$. In case of
equality, $N$ has to be $S^{\,2}_4$. So, assume $\#(V(N)) = 5$.
Define a binary relation $\sim$ on $V(N)$ by\,: $x\sim y$ if $V(N)
\setminus\{x, y\}$ is not a triangle in $N$. The hypothesis on $N$
implies that $\sim$ is an equivalence relation with at least two
equivalence classes. Since $\#(V(N)) = 5$, either there exists an
equivalence class $W$ of size 4 or $V(N)$ can be written as $V(N)
= V_1\sqcup V_2$, where $V_1$ is of size 2 and $V_1$ is a union of
equivalence classes. Accordingly, $S^{\,2}_4(W)$ or
$S^{\,0}_2(V_1)\ast S^{\,1}_3(V_2)$ is a subcomplex of $N$. \hfill
$\Box$

\begin{lemma}$\!\!${\bf .} \label{l3.3}
If a $5$-vertex simplicial complex $L$ is contractible then it is
collapsible.
\end{lemma}

\noindent {\bf Proof\,.} Let $f_i$ be the number of $i$-faces in
$L$. Suppose $\dim(L)\leq 2$. The proof is by induction on the
number $f_2$. If $f_2=0$ then, being contractible, $L$ is a tree
and hence collapsible. So assume $f_2
> 0$ and we have the result for smaller values of $f_2$. Let $N$
be the subcomplex of $L$ consisting of the triangles of $L$ and
their faces. If $N$ satisfies the hypothesis of Lemma \ref{l3.2},
then by Lemma \ref{l3.2}, ($N$ and hence) $L$ contains a
combinatorial 2-sphere as subcomplex. Then $H_2(L)\neq 0$. This is
a contradiction, since $L$ is contractible. So we may assume that
some edge of $N$ is contained in a unique triangle. Then $L$ is
collapsible to a contractible simplicial complex with fewer
triangles and hence we are done by induction hypothesis.

Consider the case when $\dim(L) = 3$. Clearly, $f_3\leq 5$. If
$f_3 = 5$ then $L = S^{\,3}_5$ and hence not contractible. So,
$f_3\leq 4$. Then there exists a triangle which is in a unique
tetrahedron. Then $L$ is collapsible to a contractible simplicial
complex with one less tetrahedron. Inductively, $L$ is collapsible
to a contractible 2-dimensional simplicial complex and hence $L$
is collapsible by the previous step.  Finally, if $\dim(L) = 4$
then $L$ consists of one 4-face and its faces. Clearly, $L$ is
collapsible in this case. \hfill $\Box$

\begin{lemma}$\!\!${\bf .} \label{l3.4}
If $M$ is a $10$-vertex combinatorial triangulation of a simply
connected $4$-manifold with $\chi(M) = 2$ then $M$ is a
combinatorial $4$-sphere.
\end{lemma}

\noindent {\bf Proof\,.} Let $\sigma$ be a facet of $M$. Let $L$
be the induced subcomplex of $M$ on $V(M)\setminus\sigma$. Then,
by Part $(a)$ of Lemma \ref{l3.1}, $L$ is contractible. Since $L$
has 5 vertices, by Lemma \ref{l3.3}, $L$ is collapsible.
Therefore, by Part $(b)$ of Lemma \ref{l3.1}, $M$ is a
combinatorial sphere. \hfill $\Box$

\bigskip

\noindent {\bf Proof of Theorem \ref{t1}\,.} Follows from Lemma
\ref{l3.4}. \hfill $\Box$

\section{Twelve-vertex complementary pseudomanifold.}

Throughout this section, $M^{\,6}_{12}$ will denote a putative
(fixed but arbitrary) 12-vertex 6-dimensional complementary weak
pseudomanifold.

\begin{lemma}$\!\!${\bf .} \label{l4.1}
$M^{\,6}_{12}$ has $12$ vertices, ${12\choose 2}=66$ edges,
${12\choose 3}= 220$ triangles, ${12\choose 4}=495$ tetrahedra,
$660$ $4$-faces, $462$ $5$-faces and $132$ facets.
\end{lemma}

\noindent {\em Proof\,.} Since $M^{\,6}_{12}$  is $6$-dimensional,
no set of $\geq 8$ vertices forms a face. Therefore, by
complementarity, $M^{\,6}_{12}$  is $4$-neighbourly. Since exactly
one set in each of the $\frac{1}{2}{12 \choose 6}$ pairs of
complementary 6-sets forms a face, it follows that the number of
5-faces is $\frac{1}{2}{12 \choose 6} =462$. Since each 5-face is
contained in two facets and each facet contains seven 5-faces, an
obvious two-way counting shows that the number of facets is
$\frac{1}{7}{12\choose 6}=132$. Finally, since a set of 5 vertices
forms a face if and only if the complementary set is not a facet,
it follows that the number of 4-faces is ${12 \choose 5}- 132 =
660$. \hfill $\Box$

\begin{defn}$\!\!${\bf .} {\rm (cf. \cite{bd3}.)}
{\rm A partition of the vertex set of $M^{\,6}_{12}$ into three
3-faces $A_1$, $A_2$, $A_3$ is called an} amicable partition {\rm
if the link of each $A_i$ is $S^{\,2}_4(A_{i+1})$ (addition in the
suffix is modulo 3)}.
\end{defn}
We have:

\begin{lemma}$\!\!${\bf .} \label{l4.2}
Let $A$ be a $3$-face of $M^{\,6}_{12}$. Suppose the link of $A$
is a standard sphere $S^{\,2}_4$. Then $A$ belongs to a unique
amicable partition of $M^{\,6}_{12}$.
\end{lemma}

\noindent {\em Proof\,.} Put $A=A_1$. Let $A_2$ be the vertex set
of the link of $A_1$ and let $A_3$ be the set of vertices outside
$A_1\cup A_2$. Then each $A_i$ contains 4 vertices. By Lemma
\ref{l4.1}, $M^{\,6}_{12}$  is $4$-neighbourly. In particular,
each $A_i$ is a 3-face of $M^{\,6}_{12}$. So, to complete the
proof, it is sufficient to show that the link of $A_2$
(respectively $A_3$) is the standard sphere on $A_3$ (respectively
$A_1$).

Take any vertex $x\in A_2$. Then $A_3\cup\{x\}$ is not a face
since its complement $A_1\cup(A_2\setminus\{x\})$ is a face. Thus
no vertex of $A_2$ belongs to the link of $A_3$. Therefore, the
vertex set of the link of $A_3$ is contained in $A_1$. Since this
link has at least 4 vertices, it follows that the link of $A_3$ is
the standard sphere on $A_1$. Replacing $A_1$ by $A_3$ (and hence
$A_2$ by $A_1$, $A_3$ by $A_2$) in this argument, we see that the
link of $A_2$ is the standard sphere on $A_3$. \hfill $\Box$

\begin{lemma}$\!\!${\bf .} \label{l4.3}
The link of any $4$-face in $M^{\,6}_{12}$ is a circle.
\end{lemma}

\noindent {\bf Proof\,.} If not, then the link (of some 4-face) is
a disconnected regular graph of degree two on at most seven
vertices. Hence the link is either $S^{\,1}_3\sqcup S^{\,1}_3$ or
$S^{\,1}_3\sqcup S^{\,1}_4$.

\medskip

\noindent {\em Case} 1. The link of a 4-face $\sigma$ is
$S^{\,1}_3(V_1) \sqcup S^{\,1}_3(V_2)$. Let $u$ be the unique
vertex outside $\sigma\cup V_1\cup V_2$. Consider
$\alpha_i=V_i\cup\{u\}$, $i=1, 2$. By Lemma \ref{l4.1}, $\alpha_i$
is a 3-face. Since $\sigma\cup e$ is a face of $M^{\,6}_{12}$ for
every 2-subset $e$ of $V_2$, it follows by complementarity that
${\rm Lk}(\alpha_1)$ has no vertex in $V_2$. Thus $V({\rm
Lk}(\alpha_1)) \subseteq\sigma$. Since ${\rm Lk}(\alpha_1)$ is a
2-dimensional weak pseudomanifold, it follows that ${\rm
Lk}(\alpha_1)$ has 4 or 5 vertices. Similarly, ${\rm
Lk}(\alpha_2)$ has 4 or 5 vertices.

\smallskip

\noindent {\em Subcase} 1.1. ${\rm Lk}(\alpha_i)$ has 4 vertices
for some $i$, say $i=1$.  By Lemma \ref{l4.2}, $\alpha_1$ belongs
to a unique amicable partition $\{\alpha_1, \beta, \gamma\}$ of
$M^{\,6}_{12}$.  Clearly, $\beta \subseteq \sigma$. This implies
that $V({\rm Lk}(\sigma))\subseteq V({\rm Lk}(\beta))=\gamma$, a
contradiction.

\smallskip

\noindent {\em Subcase} 1.2. ${\rm Lk}(\alpha_i)$ has 5 vertices
for $i= 1$, 2. Thus $V({\rm Lk}(\alpha_i))=\sigma$. Now, since the
join $S^{\,1}_3 \ast S^{\,0}_2$ of two standard spheres is the
only 5-vertex 2-dimensional weak pseudomanifold (see \cite{bd2}),
it follows that, for $i=1$, 2, ${\rm Lk}(\alpha_i) = S^{\,1}_3
\ast S^{\,0}_2$ on the common vertex set $\sigma$. Now, it is easy
to see that given any two copies of $S^{\,1}_3 \ast S^{\,0}_2$ on
a common vertex set, there is an edge of one whose complement is a
triangle in the other. Thus there must exist an edge $e$ of ${\rm
Lk}(\alpha_2)$ such that $\sigma\setminus e$ is a triangle of
${\rm Lk}(\alpha_1)$. Then the faces $\alpha_2\cup e$ and
$\alpha_1\cup(\sigma \setminus e)$ of $M^{\,6}_{12}$ cover the
vertex set, contradicting complementarity.

\medskip

\noindent {\em Case} 2. The link of a 4-face $\sigma$ is
$S^{\,1}_3 \sqcup S^{\,1}_4$. Let $\alpha$ be the vertex set of
the $S^{\,1}_4$. By Lemma \ref{l4.1}, $\alpha$ is a face of
$M^{\,6}_{12}$. Since all the edges of the $S^{\,1}_{3}$ occur in
${\rm Lk}(\sigma)$, by complementarity it follows that $V({\rm
Lk}(\alpha))\subseteq\sigma$. Therefore, ${\rm Lk}(\alpha)$ is
either an $S^{\,2}_4$ or an $S^{\,1}_3\ast S^{\,0}_2$.

\smallskip

\noindent {\em Subcase} 2.1. ${\rm Lk}(\alpha)=S^{\,2}_4$. Then,
by Lemma \ref{l4.2}, $\alpha$ belongs to a unique amicable
partition $\{\alpha, \beta, \gamma\}$ of $M^{\,6}_{12}$.  Clearly,
$\beta\subseteq\sigma$ and hence $V({\rm Lk}(\sigma))\subseteq
V({\rm Lk}(\beta)) = \gamma$, contradiction.

\smallskip

\noindent {\em Subcase} 2.2. ${\rm Lk}(\alpha)=S^{\,1}_3\ast
S^{\,0}_2$. Thus $V({\rm Lk}(\alpha))=\sigma$. Let ${\rm
Lk}(\sigma)=S^{\,1}_3(V_1) \sqcup(S^{\,0}_2(\{x, y\})\ast
S^{\,0}_2(\{z, w\}))$. Choose one of the four vertices $x$, $y$,
$z$, $w$, say we choose $x$. Consider the 3-face $\beta=
V_1\cup\{x\}$. Since $yz$ and $yw$ are edges in ${\rm
Lk}(\sigma)$, complementarity implies that $z$ and $w$ are not
vertices of ${\rm Lk}(\beta)$. So, $V({\rm
Lk}(\beta))\subseteq\sigma\cup\{y\}$. Thus ${\rm Lk}(\beta)$ is a
2-dimensional weak pseudomanifold on $\leq 6$ vertices. Let $\{a,
b\}\subseteq\sigma$ be the unique non-edge in ${\rm Lk}(\alpha)$.
Since ${\rm Lk}(\alpha)$ contains three triangles through $a$,
complementarity implies that at least three edges through $b$
(contained in $\si$) are missing (in ${\rm Lk}(V_1)$ and hence) in
${\rm Lk}(\beta)$. Since any vertex in a 2-dimensional weak
pseudomanifold on $\leq 6$ vertices can be on at most two
non-edges, it follows that $b$ is not a vertex of ${\rm
Lk}(\beta)$. Similarly, $a$ is not a vertex of ${\rm Lk}(\beta)$.
Therefore $V({\rm Lk}(\beta))\subseteq\si\cup\{y\}\setminus\{a,
b\}$ and hence ${\rm
Lk}(\beta)=S^{\,2}_4(\si\cup\{y\}\setminus\{a, b\})$. In
particular, $\si\setminus\{a, b\}$ is a face of ${\rm Lk}(\beta)$.
Hence $V_1\cup(\si\setminus\{a, b\})\cup\{x\}$ is a facet of
$M^{\,6}_{12}$. Since this argument goes through with any of the
four vertices of $\alpha$ in place of $x$, this shows that the
5-face $V_1 \cup\si\setminus\{a, b\}$ is contained in (at least)
four facets of $M^{\,6}_{12}$. This is a contradiction since
$M^{\,6}_{12}$ is a weak pseudomanifold. \hfill $\Box$

\begin{lemma}$\!\!${\bf .} \label{l4.4}
The link of any tetrahedron in $M^{\,6}_{12}$ is a connected
combinatorial $2$-manifold.
\end{lemma}

\noindent {\bf Proof\,.} Let $L$ be the link of a 3-face $\si$.
Then $L$ is a 2-dimensional weak pseudomanifold on at most 8
vertices. By Lemma \ref{l4.3}, the link in $L$ of each vertex is a
circle. Hence $L$ is a combinatorial 2-manifold. If $L$ were
disconnected, (since any 2-dimensional weak pseudomanifold has
$\geq 4$ vertices with equality only for $S^{\,2}_4$) $L$ would
have to be $S^{\,2}_4\sqcup S^{\,2}_4$. Say,
$L=S^{\,2}_4(\alpha)\sqcup S^{\,2}_4(\beta)$. By complementarity,
${\rm Lk}(\alpha)=S^{\,2}_4(\sigma)={\rm Lk}(\beta)$. This
contradicts Lemma \ref{l4.2}. \hfill $\Box$

\begin{lemma}$\!\!${\bf .} \label{l4.5}
The link of any tetrahedron in $M^{\,6}_{12}$ is a combinatorial
$2$-sphere.
\end{lemma}

\noindent {\bf Proof\,.} For $i\geq 4$, let $c_i$ be the number of
3-faces of degree $i$ in $M^{\,6}_{12}$. Counting the total number
of 3-faces in $M^{\,6}_{12}$ and the total number of pairs
$(\alpha, \beta)$ where $\alpha\subseteq\beta$, $\alpha$ is a
3-face and $\beta$ is a 4-face of $M^{\,6}_{12}$ we get (in view
of Lemma \ref{l4.1})\,:
\begin{equation}
\sum c_i  =  495 ~~ \mbox{ and } ~~ \sum i c_i  = 660\times 5 =
3300. \label{eq1}
\end{equation}
Now, for any 3-face $\alpha$ of $M^{\,6}_{12}$ , the link of
$\alpha$ is a connected combinatorial 2-manifold (by Lemma
\ref{l4.4}) and hence has Euler characteristic $\leq 2$. If this
link has $i$ vertices then it follows that it has $\geq 2i-4$
triangles, i.e., $\alpha$ is contained in at least $2i-4$ facets.
Hence $\sum(2i-4) c_i$ is a lower bound on the number of pairs
$(\alpha, \gamma)$, where $\alpha\subseteq\gamma$, $\alpha$ is a
3-face and $\gamma$ is a facet of $M^{\,6}_{12}$. But by Lemma
\ref{l4.1}, the number of such pairs is $132\times 35$. Therefore
we get
$$
\sum(2i - 4) c_i \leq 132 \times 35 = 4620.
$$
By (\ref{eq1}), equality holds in this inequality. Therefore,
equality holds throughout the above argument. Thus the link of
each 3-face is a combinatorial 2-manifold of Euler characteristic
2, and hence it is a combinatorial 2-sphere. \hfill $\Box$

\begin{lemma}$\!\!${\bf .} \label{l4.6}
The link of any triangle in $M^{\,6}_{12}$ is a connected
combinatorial $3$-manifold on $9$ vertices.
\end{lemma}

\noindent {\bf Proof\,.} The link is a combinatorial 3-manifold by
Lemma \ref{l4.5}. It has 9 vertices by Lemma \ref{l4.1}. It is
connected since any disconnected combinatorial 3-manifold needs at
least $5+5=10$ vertices. \hfill $\Box$

\begin{lemma}$\!\!${\bf .} \label{l4.7}
If $\delta$ is a triangle in $M^{\,6}_{12}$ such that $L = {\rm
Lk}(\delta)$ is $2$-neighbourly then $L$ does not have any induced
$S^{\,2}_4$.
\end{lemma}

\noindent {\bf Proof\,.} If possible let $S^{\,2}_4(X)$ be an
induced sub-complex of $L$. Put $Y=V(L)\setminus X$. By
complementarity, $Y$ is a 4-face of $M^{\,6}_{12}$ and ${\rm
Lk}(Y)$ does not contain any vertex of $X$. This implies that
${\rm Lk}(Y)=S^{\,1}_3(\delta)$. Therefore, for each vertex $x$ in
$\delta$, $Y\cup\delta\setminus\{x\}$ is a facet of
$M^{\,6}_{12}$. Hence, by complementarity, $x\not\in {\rm Lk}(X)$.
Since this holds for each $x\in\delta$, $\delta\cap V({\rm
Lk}(X))=\emptyset$. So, if $\tau$ is a facet of $M^{\,6}_{12}$
containing $X$ then $\tau$ is disjoint from $\delta$. Thus
$\alpha=V(M^{\,6}_{12})\setminus \tau$ is a set of 5 vertices
containing $\delta$ and, by complementarity, $\alpha$ is not a
face. Thus $\alpha\setminus\delta$ is not an edge of ${\rm
Lk}(\delta)$. This contradicts the assumption that ${\rm
Lk}(\delta)$ is 2-neighbourly. \hfill $\Box$

\begin{lemma}$\!\!${\bf .} \label{l4.8}
The link of any triangle in $M^{\,6}_{12}$ is a combinatorial
$3$-sphere on $9$ vertices.
\end{lemma}

\noindent {\bf Proof\,.} By Proposition \ref{alst}, there is a
unique 9-vertex combinatorial 3-manifold which is not a
combinatorial 3-sphere, namely the 3-dimensional Klein bottle
$K^{\,3}_9$. So, in view of Lemma \ref{l4.6}, it is enough to
prove that $K^{\,3}_9$ can not occur as the link of a triangle of
$M^{\,6}_{12}$. Now, $K^{\,3}_9$ is 2-neighbourly and has an
induced subcomplex isomorphic to $S^{\,2}_4$. Indeed any $4$-path
in the underlying $S^{\,1}_9$ induces an $S^{\,2}_4$ in
$K^{\,3}_9$. The result now follows from Lemma \ref{l4.7}. \hfill
$\Box$

\begin{lemma}$\!\!${\bf .} \label{l4.9}
$M^{\,6}_{12}$ is a pseudomanifold.
\end{lemma}

\noindent {\bf Proof\,.} By complementarity, two facets of
$M^{\,6}_{12}$ can not cover its vertex set. Hence any two of the
facets have at least a triangle in common. Since by Lemma
\ref{l4.6}, the link of any triangle is strongly connected, it
follows that $M^{\,6}_{12}$ is strongly connected and hence is a
pseudomanifold. \hfill $\Box$

\begin{lemma}$\!\!${\bf .} \label{l4.10}
Let $\alpha$ be a facet of $M^{\,6}_{12}$. For $0\leq j\leq 2$ let
$e_j$ be the number of $(6-j)$-faces meeting $\alpha$ in exactly
$6-j$ vertices. Then $e_0=7$, $e_1=51$ and $e_2=139$.
\end{lemma}

\noindent {\bf Proof\,.} Let ${\cal A}$ be the set of all
$3$-faces of $M^{\,6}_{12}$ disjoint from $\alpha$. For $\gamma
\in{\cal A}$, let $d_j(\gamma)$ be the number of $j$-faces in
${\rm Lk}(\gamma)$, $0\leq j\leq 2$. Clearly $\sum_{\gamma\in{\cal
A}} d_j(\gamma)$ counts the number of pairs $(\gamma, \delta)$
where $\gamma\in{\cal A}$ and $\delta$ is a $(4+j)$-face of
$M^{\,6}_{12}$ containing $\gamma$. Since $\alpha$ is a facet of
$M^{\,6}_{12}$, by complementarity, $V({\rm Lk}(\gamma))
\subseteq\alpha$ for $\gamma\in{\cal A}$. It follows that $\gamma=
\delta\setminus\alpha$ for any such pair $(\delta, \gamma)$.
Therefore, the number of such pairs $(\gamma, \delta)$ equals the
number of $(4+j)$-faces $\delta$ of $M^{\,6}_{12}$ meeting
$\alpha$ in $j+1$ vertices. By complementarity, this equals the
number of vertex sets of size $7-j$ meeting $\alpha$ in $6-j$
vertices, which do not form a face of $M^{\,6}_{12}$. Out of the
$5\times {7\choose j+1}$ sets of size $7-j$ meeting $\alpha$ in
$6-j$ vertices, exactly $e_j$ are faces of $M^{\,6}_{12}$. So we
get
\begin{equation}
e_j=5\times {7\choose j+1}-\sum_{\gamma\in{\cal A}} d_j(\gamma),
~~~ 0\leq j\leq 2. \label{eq2}
\end{equation}
Now, for each fixed $\gamma\in{\cal A}$, $d_j(\gamma)$ is the
numbers of $j$-faces in ${\rm Lk}(\gamma)$. Since ${\rm
Lk}(\gamma)$ is an $S^{\,2}$ (Lemma \ref{l4.5}), we have
\begin{equation}
d_1(\gamma) = 3(d_0(\gamma)-2), ~~ d_2(\gamma) = 2(d_0(\gamma)-2).
\label{eq3}
\end{equation}
Adding (\ref{eq3}) over all $\gamma\in{\cal A}$, and plugging the
result into (\ref{eq2}) we get $e_1=3e_0+30$ and $e_2=125+2e_0$.
Since  $M^{\,6}_{12}$ is a $6$-dimensional weak pseudomanifold, we
have $e_0=7$ and hence $e_1=51$ and $e_2=139$. \hfill $\Box$

\begin{lemma}$\!\!${\bf .} \label{l4.11}
Given any facet of $M^{\,6}_{12}$ the number of facets meeting the
given facet in $3$, $4$, $5$, $6$ vertices equals $36$, $58$, $30$
and $7$ respectively.
\end{lemma}

\noindent {\bf Proof\,.} Fix a facet $\alpha$ of $M^{\,6}_{12}$.
Clearly exactly $7$ facets meet $\alpha$ in $6$ vertices.

By Lemma \ref{l4.10} (with $j=2$), exactly $139$ \, $4$-faces of
$M^{\,6}_{12}$ meet $\alpha$ in $4$ vertices.  Therefore, by
complementarity, out of the $5\times {7\choose 4}$ sets of size 7
meeting $\alpha$ in 3 vertices, exactly 139 are not facets. So,
$5\times{7\choose 4} - 139 = 36$ facets meet $\alpha$ in a
triangle.

Let ${\cal B}$ denote the set of all $4$-faces contained in
$\alpha$. For $\gamma\in {\cal B}$ and $0\leq j\leq 1$, let
$d_j(\gamma)$ be the number of $j$-faces in ${\rm Lk}(\gamma)$.
Lemma \ref{l4.10} (with $j=1$) shows that $\sum_{\gamma\in {\cal
B}} (d_0(\gamma) - 2) = 51$. That is, $\sum_{\gamma\in {\cal B}}
d_0(\gamma) = 51 + 2 \times{7\choose 5} = 93$.  Since ${\rm
Lk}(\gamma)$ is a circle, $d_0(\gamma) = d_1(\gamma)$.  Hence
$\sum_{\gamma\in {\cal B}} d_1(\gamma) = 93$, and therefore
$\sum_{\gamma\in {\cal B}} (d_1(\gamma)- 3) = 93 -
3\times{7\choose 5} = 30$.  But this last sum counts the number of
ordered pairs $(\gamma , e)$ where $\gamma \in {\cal B}$ and $e$
is an edge of ${\rm Lk}(\gamma)$ disjoint from $\alpha$. Since the
number of such pairs equals the number of facets $\beta =
\gamma\cup e$ meeting $\alpha$ in exactly $5$ vertices, we see
that there are $30$ such $\beta$'s.

Since, by Lemma \ref{l4.1}, there are $131$ facets other than
$\alpha$, we find by subtraction that $131-36-30-7 = 58$ facets
meet $\alpha$ in $4$ vertices. \hfill $\Box$

\begin{lemma}$\!\!${\bf .} \label{l4.12}
$(a)$ Each edge   of $M^{\,6}_{12}$ is contained in exactly $42$
facets.
\newline \mbox{} \hspace{24.5mm}
$(b)$ Each vertex of $M^{\,6}_{12}$ is contained in exactly $77$
facets.
\end{lemma}

\noindent {\bf Proof\,.} For $i\geq 0$, let $a_i$ be the number of
edges of $M^{\,6}_{12}$ which are contained in exactly $i$ facets.
Counting in two ways (i) total number of edges, (ii) total number
of pairs $(e, \alpha)$ where $e$ is an edge, $\alpha$ is a facet
and $e\subseteq \alpha$ and (iii) total number of triples $(e,
\alpha_1, \alpha_2)$ where $e$ is an edge, $\alpha_1$, $\alpha_2$
are distinct facets and $e\subseteq\alpha_1\cap\alpha_2$, we get
(in view of Lemmas \ref{l4.1} and \ref{l4.11})\,:
\begin{eqnarray*}
\sum a_i & = & {12\choose 2} = 66, \hspace{10mm}
\sum i\,a_i   =   132 \times {7\choose 2} = 2772, \\
\sum i(i-1)\,a_i & = & 132 \times \left(36 \times {3\choose 2} +
58 \times {4\choose 2} + 30 \times {5\choose 2} + 7 \times
{6\choose 2}\right) =  113652.
\end{eqnarray*}
Hence we find
$$\sum (i - 42)^2 a_i = 0.
$$
That is, $a_i = 0$ for $i\neq 42$. So every edge of $M^{\,6}_{12}$
is contained in exactly $42$ facets.

Next, fix a vertex $x$ of $M^{\,6}_{12}$. Counting in two ways the
number of pairs $(e, \alpha)$ where $e$ is an edge containing $x$
and $\alpha$ is a facet containing $e$, we find that the number
$r$ of facets through $x$ satisfies $r \times 6 = 11 \times 42$.
So, $r=77$. Thus each vertex is in $77$ facets. \hfill $\Box$

\begin{lemma}$\!\!${\bf .} \label{l4.13}
The link of any edge in $M^{\,6}_{12}$ is a $10$-vertex simply
connected combinatorial $4$-manifold of Euler characteristic $2$.
\end{lemma}

\noindent {\bf Proof\,.} Let $L$ be the link of an edge. By Lemma
\ref{l4.1}, $L$ is a $10$-vertex $2$-neighbourly weak
pseudomanifold. By Lemma \ref{l4.8}, $L$ is a combinatorial
$4$-manifold.  Being 2-neighbourly, $L$ is connected. By
Proposition \ref{brku2}, $|L|$ is simply connected.

If $e$ is an edge of $M^{\,6}_{12}$ then, by Lemma \ref{l4.12},
$77\times 2 - 42 = 112$ facets intersect $e$ and hence $132 - 112
= 20$ facets are disjoint from $e$.  Therefore, by
complementarity, each edge of $M^{\,6}_{12}$ is contained in
${10\choose 3} - 20 = 100$ \, $4$-faces.

For $0\leq i\leq 4$, let $f_i$ be the number of $i$-faces of $L$.
Since $L$ is $2$-neighbourly with $10$ vertices, we have $f_0 =
10$ and $f_1={10\choose 2} = 45$. By Lemma \ref{l4.12}, $f_4=42$
and, by the above argument, $f_2 = 100$. An obvious two way
counting yields $f_3 = 5f_4/2 = 105$.  Hence the Euler
characteristic of $L$ is $10-45+100-105+ 42= 2$. \hfill $\Box$

\begin{lemma}$\!\!${\bf .} \label{l4.14}
$M^{\,6}_{12}$ is a combinatorial manifold.
\end{lemma}

\noindent {\bf Proof\,.} By Lemmas \ref{l4.13} and \ref{l3.4}, the
link of each edge in $M^{\,6}_{12}$ is a combinatorial 4-sphere.
Then, by Lemma \ref{l4.1}, the link of any vertex is an 11-vertex
combinatorial 5-manifold. Now Part (a) of Proposition \ref{brku1}
implies that the link of any vertex is a combinatorial 5-sphere.
Hence $M^{\,6}_{12}$ is a combinatorial 6-manifold. \hfill $\Box$

\bigskip

\noindent {\bf Proof of Theorem \ref{t2}\,.}  By Remark \ref{re1}
(or by Lemma \ref{l4.1}) $\chi(M^{\,6}_{12}) = 1$. Since
$M^{\,6}_{12}$ is 3-neighbourly, $|M^{\,6}_{12}|$ is simply
connected and hence orientable. This contradicts Lemma \ref{l4.14}
since the Euler characteristic of an orientable closed manifold of
dimension $\equiv 2$ (mod 4) is even (see \cite[Corollary
26.11]{gh}). (Lemma \ref{l4.14} also contradicts Proposition
\ref{da1} and Proposition \ref{brku1} $(b)$.) \hfill $\Box$

\bigskip

\noindent {\bf Acknowledgement\,:} Lemma \ref{l3.1} is essentially
due to W. K\"{u}hnel. The authors thank W. K\"{u}hnel for pointing
it out and for suggesting the validity of Lemma \ref{l3.3} in an
e-mail message to the second author. The authors are thankful to
S. Gadgil and V. Pati for numerous useful conversations. The
authors thank the anonymous referee for suggesting the inclusion
of Lemma \ref{l3.4} so that the results can be proved without
using the 4-dimensional Poincar\'{e} conjecture, as well as for
many other useful comments which helped to improve the
presentation of this paper. Remark \ref{re2} is also due to the
referee.

\bigskip

{\footnotesize

\noindent B. Bagchi, Statistics-Mathematics Unit, Indian
Statistical Institute, Bangalore 560\,059, India.
\newline \mbox{} \hspace{5mm} E-mail: bbagchi@isibang.ac.in

\noindent B. Datta, Department of Mathematics, Indian Institute of
Science, Bangalore 560\,012, India. \newline \mbox{} \hspace{5mm}
E-mail: dattab@math.iisc.ernet.in }

\end{document}